\documentstyle[12pt]{article}

\def\Hom{\mathop{\rm Hom}\nolimits}

\def\Ad{\mathop{\rm Ad}\nolimits}
\def\Ext{\mathop{\rm Ext}\nolimits}
\def\id{\mathop{\rm id}\nolimits}
\def\Bott{\mathop{\rm Bott}\nolimits}

   \newcommand{\ah}{asymptotic homomorphism }
   \newtheorem{thm}{Theorem}[section]
   \newtheorem{prop}[thm]{Proposition}
   \newtheorem{lemma}[thm]{Lemma}
   \newtheorem{cor}[thm]{Corollary}

\title{Asymptotically split extensions and $E$-theory}
\author{Vladimir Manuilov${}^1$ and Klaus Thomsen}

\date{}

\begin{document}

\maketitle\footnotetext[1]{%
Partially supported by RFFI, grant No 99-01-01201}

\begin{abstract}
We show that the $E$-theory of Connes and Higson can be
formulated in terms of $C^*$-extensions in a way quite similar to the way
in which the $KK$-theory of Kasparov can. The essential difference is that
the role played by split extensions should be taken by asymptotically split
extensions. We call an extension of a $C^*$-algebra $A$ by a stable
$C^*$-algebra $B$ {\it asymptotically split} if there exists an asymptotic
homomorphism consisting of right inverses for the quotient map.
An extension is called {\it semi-invertible} if it can be made
asymptotically split by adding another extension to it. Our main result is
that there exists a one-to-one correspondence between asymptotic
homomorphisms from $SA$ to $B$ and homotopy classes of semi-invertible
extensions of $S^2A$ by $B$.

\end{abstract}

\section{Introduction}

Connes and Higson introduced in \cite{CH} a
construction which produces an asymptotic homomorphism out of an extension
of $C^*$-algebras. The Connes-Higson construction is the backbone of
$E$-theory and gives us a way to study $C^*$-extensions via asymptotic
homomorphisms.
Such a translation can be quite powerful within the
territory of $KK$-theory, where the $C^*$-extensions are semi-split, i.e.
admit a completely positive contraction as a right-inverse for the quotient
map. It is namely known that the Connes--Higson construction sets up a
bijection between homotopy classes of semi-split extensions and completely
positive asymptotic homomorphisms. This bijection is particularly useful
because completely positive asymptotic homomorphisms are easier to handle
than general ones, and because the powerful homotopy invariance results of
Kasparov, \cite{K}, allows one to translate homotopy information to more
algebraic information about the $C^*$-extensions. This well-behaved
correspondance between semi-split $C^*$-extensions and homotopy classes of
completely positive asymptotic homomorphisms was used in \cite{MT} to
obtain a better understanding of the short exact sequence of the
UCT-theorem by identifying the kernel of the map from $KK(A,B) =
\Ext^{-1}(SA,B)$ to $KL(A,B)$ as the group arising from the weakly
quasi-diagonal extensions of $SA$ by $B \otimes {\cal K}$.

The present paper originated in the
desire to extend the nice relation between $C^*$-extensions and asymptotic
homomorphisms beyond the case of semi-split extensions. The key problem in
this connection is (at least for the moment) to decide if the Connes--Higson
construction is injective in general. In other words, the problem is to
decide if two $C^*$-extensions - with stable and maybe suspended ideals -
which give rise to homotopic asymptotic homomorphisms must themselves be
homotopic. From \cite{HL-T} we know that this is the case when both
extensions are suspensions and the result of the present paper shows
that it is also the case when both extensions are what we call
semi-invertible and the quotient $C^*$-algebra is a double extension. But
in general we still don't know the answer. Nonetheless, we shall show here
that there is a way to faithfully represent $E$-theory by use of
$C^*$-extensions which does not require infinitely many suspensions as in
\cite{HL-T} or longer decomposition series as in \cite{Cu}.

To describe this, let $A$ and $B$ be separable
$C^*$-algebras and assume for simplicity that $B$ is stable. We call an
extension of $A$ by $B$ {\it asymptotically split} when there is a family
$(\pi_t)_{t\in[1,\infty)}$ of right-inverses for the quotient map such
that $(\pi_t)_{t\in[1,\infty)}$ is an asymptotic homomorphism. An
extension is then {\it semi-invertible} when it can be made asymptotically
split by adding another extenson to it. We prove that

\vspace{-1\itemsep}
\begin{enumerate}
\item[1)] Every asymptotic homomorphism $S^2 A \to B$ is homotopic to the
Connes--Higson construction of a semi-invertible extension of $SA$ by $B$.

\vspace{-1.2\itemsep}
\item[2)] Two semi-invertible extensions of $S^2A$ by $B$ are homotopic
(as semi-invertible extensions) if and only if the Connes--Higson
construction applied to them give homotopic asymptotic homomorphisms.
\end{enumerate}

\vspace{-1.2\itemsep}

These results show that the $E$-theory of Connes and Higson can be
formulated in terms of $C^*$-extensions in a way quite similar to the
way in which the $KK$-theory of Kasparov can. The essential difference
is that the role played by split extensions should be taken by
asymptotically split extensions. It is our hope that this parallel
between the way $KK$-theory and $E$-theory can be described in terms
of $C^*$-extensions can be strenghtened even further. In particular
it would be nice if some of the suspensions occuring in 1) and 2) could be
removed and if one could substitute homotopy with a more algebraic
relation in the description of $E$-theory.

\section{Asymptotically split extensions and $\Ext^{- 1/2}$}

In the following $A$ and $B$ are separable $C^*$-algebras, $B$ stable,
i.e. $B=B\otimes\cal K$, where $\cal K$ denotes the $C^*$-algebra
of compact operators. As usual, we denote by $C_0(X)$ the
$C^*$-algebra of continuous functions on $X$ vanishing at infinity,
and $SA=C_0(0,1)\otimes A$ denotes the suspension $C^*$-algebra
over $A$.
Let $M(B)$ denote the multiplier algebra of $B$, \cite{P},
$Q(B) = M(B)/B$ the
corresponding corona algebra and $q_B : M(B) \to Q(B)$ the quotient map.
We shall identify
the set of extensions of $A$ by $B$ with $\Hom(A,Q(B))$ in the standard way, \cite{B}.
Two extensions $\varphi, \ \psi : A \to Q(B)$ are {\it unitarily equivalent}
when there is a unitary $w \in M(B)$ such that
$\Ad q_B(w) \circ \varphi = \psi$. As is wellknown the set of unitary
equivalence classes of extensions of $A$ by $B$ form a semi-group and
we denote this semi-group by $\Ext(A,B)$.

Recall that an asymptotic homomorphism from $A$ to
$B$ is a family $\varphi=\{\varphi_t\}_{t\in [1,\infty)}:A\to B$ of
maps such that $t \mapsto \varphi_t(a)$ is continuous for any $a\in A$ and the $\varphi_t$'s behave like a $*$-homomorphism
asymptotically as $t\to\infty$, \cite{CH}.
Namely, for any $a,b\in A$,
$\lambda\in\bf C$ one has
$$
\lim_{t \to \infty} \|\varphi_t(a^*) - \varphi_t(a)^*\| = 0,
$$
$$
\lim_{t\to\infty}\|\varphi_t(\lambda a+b)-\lambda\varphi_t(a)-
\varphi_t(b)\|=0,
 $$
 $$
\lim_{t\to\infty}\|\varphi_t(ab)-\varphi_t(a)\varphi_t(b)\|=0.
 $$
Two asymptotic homomorphisms $\varphi$ and $\psi$ are equivalent when
$\lim_{t\to\infty}\|\varphi_t(a)-\psi_t(a)\|=0$ for all $a\in A$ and
are homotopic when  there exists an asymptotic homomorphism
$\phi=\{\phi_t\}_{t\in[1,\infty)}: A\to C[0,1]\otimes B$ such that
the compositions with the evaluation maps at $0$ and $1$ coincide with
$\varphi$ and $\psi$, respectively. The semi-group of homotopy
classes of asymptotic homomorphisms we denote by $[[A,B]]$.

An extension
$\varphi : A  \to Q(B)$ is called {\it asymptotically split}
when there is an asymptotic homomorphism
$\pi = \{\pi_t\}_{t \in [1,\infty)} : A \to M(B)$ such that
$q_B \circ \pi_t = \varphi$ for all $t$. An extension
$\varphi : A \to Q(B)$ is called {\it semi-invertible} when there is an
extension $\psi$ such that $\varphi \oplus \psi : A \to Q(B)$ is
asymptotically split. Two semi-invertible extensions are called
{\it stably equivalent} when they become unitarily equivalent after
addition by asymptotically split extensions.

Stable equivalence is an equivalence
relation on the subset of semi-invertible extensions in $\Hom(A,Q(B))$
and the corresponding equivalence classes form an abelian group which
we denote by $\Ext^{ -1/2}(A,B)$. $\Ext^{- 1/2}$ is a bifunctor which
is contravariant in the first variable, $A$, and covariant with respect
to quasi-unital $*$-homomorphisms in the second variable, $B$. It is
easy to see that the Connes--Higson construction, \cite{CH}, annihilates
asymptotically split extensions and therefore gives rise to a group
homomorphism
$$
CH : \Ext^{ -1/2} (A,B) \to [[SA,B]] \ .
$$
Two semi-invertible extensions
$$
0 {\ \longrightarrow\ }    B {\ \longrightarrow\ }
E_1 {\ \longrightarrow\ }    A {\ \longrightarrow\ }    0
$$
and
$$
0 {\ \longrightarrow\ }    B {\ \longrightarrow\ }
E_2 {\ \longrightarrow\ }    A {\ \longrightarrow\ }    0
$$
are called {\it homotopic} when there is a commuting diagram
$$
\begin{array}{ccccccccc}
0& {\ \longrightarrow\ } & B & {\ \longrightarrow\ } &
E_1 & {\ \longrightarrow\ } & A & {\ \longrightarrow\ } & 0 \\
&&{\scriptstyle \pi_0}\uparrow{\phantom a} && \uparrow &&\Vert&&\\
0& {\ \longrightarrow\ }  & \!\!\!\! C[0,1] \otimes B \!\!\!\! &
{\ \longrightarrow\ } & E & {\ \longrightarrow\ } & A &
{\ \longrightarrow\ } & 0\\
&&{\scriptstyle \pi_1}\downarrow{\phantom a} && \downarrow&&\Vert&&\\
0& {\ \longrightarrow\ } & B & {\ \longrightarrow\ } &
E_2 & {\ \longrightarrow\ } & A & {\ \longrightarrow\ } & 0
\end{array}
$$
of semi-invertible extensions. The $*$-homomorphisms
$\pi_0, \ \pi_1 : C[0,1] \otimes B \to B$ are here the surjections
obtained from evalution at the endpoints of $[0,1]$.

\bigskip

The main tool in this paper is the map $E$ introduced in
\cite{MM}, cf. \cite{MT}.
We recall the construction here. Given an \ah
$\varphi = \{\varphi_t\}_{t \in [1,\infty)} : A \to B$ we choose a sequence
$1 \leq t_1 \leq t_2 \leq t_3 \leq \cdots $ such that
$$
\lim_{i \to \infty} t_i = \infty
\quad {\rm and}\quad
\lim_{i \to \infty}
\sup_{t  \in [t_i,t_{i+1}]} \|\varphi_t(a) - \varphi_{t_i}(a)\| = 0
$$
for all $a \in A$. Let $e_{ij}, i,j \in \bf Z$ denote
the standard matrix units, which act on the standard
Hilbert $B$-module $l_2({\bf Z})\otimes B$ in the obvious way. Then
$$
\Phi (a) = \sum_{i \geq 1} \varphi_{t_i}(a) e_{ii}
$$
defines a map $\Phi : A \to {\cal L}_B(l_2({\bf Z}) \otimes B)$, where
${\cal L}_B(l_2({\bf Z}) \otimes B)$ is the $C^*$-algebra of bounded
adjointable operators on the Hilbert $C^*$-module $l_2({\bf Z}) \otimes B$.
We identify ${\cal K} \otimes B$ with the ideal of $B$-compact operators in
${\cal L}_B(l_2({\bf Z}) \otimes B)$ and observe that $\Phi$ is a
$*$-homomorphism modulo ${\cal K} \otimes B$. Furthermore,
$\Phi(a)$ commutes modulo ${\cal K} \otimes B$ with the two-sided
shift $T = \sum_{j \in {\bf Z}} e_{j,j+1}$. So we get in this way a
$*$-homomorphism
$$
E(\varphi) : C({\bf T}) \otimes A \to Q({\cal K} \otimes B) =
{\cal L}_B(l_2({\bf Z})
\otimes B)/{\cal K} \otimes B \
$$
such that
$$
E(\varphi)(f \otimes a) = f(\underline{T})\underline{\Phi(a)} \ ,
\ f \in C({\bf T}), \ a \in A \ .
$$
Here and in the following we denote by $\underline{S}$ the image in
$Q({\cal K} \otimes B) =  {\cal L}_B(l_2({\bf Z}) \otimes B)/
{\cal K} \otimes B$ of an element $S \in {\cal L}_B(l_2({\bf Z})
\otimes B)$.
It can be checked directly that the map $E$ is well-defined and does not
depend on the choice of a discretization.

\begin{lemma}\label{L1}
$E(\varphi) \in \Ext^{ - 1/2}(C({\bf T})\otimes A,
{\cal K} \otimes B)$.
\end{lemma}

{\bf Proof.} Let $-E(\varphi) : C({\bf T}) \otimes A \to Q({\cal K}
\otimes B)$ be the extension which results when we in the construction
of $E(\varphi)$ use
$$
\Psi(a) = \sum_{i \leq 0} \varphi_{t_{ - i +1}}(a)e_{ii}
$$
instead of $\Phi$. Then $-E(\varphi) \oplus E(\varphi)$ is unitary
equivalent to an extension $\psi : C({\bf T}) \otimes A \to Q({\cal K}
\otimes B)$ such that
$\psi(f\otimes a) = \underline{\pi_t(f\otimes a)}$ for all $t \in
[1,\infty)$, $f\in C({\bf T})$, $a\in A$,
where
$$
\pi=\{\pi_t\}_{t \in [1,\infty)} : C({\bf T}) \otimes A \to {\cal
L}_B(l_2({\bf Z}) \otimes B)
$$
is an \ah
obtained by convex interpolation of maps $\pi_n, \ n \in \bf N$, with
the property that
$$
 \pi_n( f \otimes a) - f(T)
\Bigl(\sum_{|i| \leq n} \varphi_{t_n}(a)e_{ii} +
\sum_{i > n} \varphi_{t_i}(a)e_{ii} + \sum_{ i < -n}
\varphi_{t_{- i +1}}(a)e_{ii}\Bigr)  \in {\cal K} \otimes B
$$
and
$$
\lim_{ n \to \infty} \pi_n( f \otimes a) - f(T)
\Bigl(\sum_{|i| \leq n}
\varphi_{t_n}(a)e_{ii} + \sum_{i > n} \varphi_{t_i}(a)e_{ii} +
\sum_{ i < -n} \varphi_{t_{- i +1}}(a)e_{ii}\Bigr) = 0 \ ,
$$
$ f \in C({\bf T}), \ a \in A$. But $\psi$ is obviously asymptotically
split.
\hfill $\Box$

\medskip

Let $\Ext^{-1/2}(A,B)_h$ denote the abelian group of homotopy classes of
semi-invertible extensions of $A$ by $B$. $\Ext^{-1/2}(A,B)_h$ is then a
quotient of $\Ext^{-1/2}(A,B)$. By homotopy invariance of the Connes--Higson
construction we get a map
$$
CH : \Ext^{-1/2}(A,B)_h \to [[SA,B]].
$$
Thanks to Lemma \ref{L1} we get from the above construction a well-defined
map
$$
E : [[A,B]] \to \Ext^{-1/2}(C({\bf T}) \otimes A,B)_h \ ,
$$
cf. \cite{MT}. By pulling back along the canonical inclusion
$S A \subseteq C({\bf T}) \otimes A$ we can also consider $E$ as a map
$$
E : [[A,B]] \to \Ext^{-1/2}(SA,B)_h.
$$
Our main result can now be
formulated as follows.

\begin{thm}\label{MAIN}

\bigskip

a$)$ The map $CH : \Ext^{-1/2}(SA,B) \to [[S^2A,B]]$ is surjective.

\medskip

b$)$ The map $E : [[SA,B]] \to \Ext^{-1/2}(S^2A,B)_h$ is an isomorphism.

\end{thm}

Let $\chi : [[SA,B]] \to [[S^2C(\bf T)\otimes A,B]]$ be the map obtained
by taking the exterior product product with
the \ah  $S^2C(\bf T) \to S\otimes {\cal K}$ which is the suspension of the
\ah $S^2\to\cal K$ obtained by applying the Connes--Higson construction
to the Toeplits extension. The composition of $\chi$ with the obvious map
$[[S^2C(\bf T)\otimes A,B]] \to [[S^3A,B]]$ will be denoted by $\chi_0$.
To prove a) we use the following statement, cf.\cite{MT,MM}.
\begin{lemma}
The diagram
\begin{equation}\label{Diag1}
\begin{array}{ccc}
\Ext^{-1/2} (S^2A,B)_h\!\!\!\!\!\!\!\!\!\!\!\!\!\!\!\!\! & & \\
{\scriptstyle E}\uparrow\phantom{a} &
\phantom{aaaa}\searrow\!{\scriptstyle CH} & \\
{[[SA,B]]} & \!\!\!\!\!\!\!\stackrel{\chi}{\longrightarrow} &
\!\!\!\!\!\!\! {[[S^3A,B]]} \, ,
\end{array}
\end{equation}
commutes.
\end{lemma}

{\bf Proof.}
We are going to prove
commutativity of the diagram
 $$
\begin{array}{ccc}
\Ext^{-1/2} (C({\bf T})\otimes A,B)_h
\!\!\!\!\!\!\!\!\!\!\!\!\!\!\!\!\!\!\!\!\!\!\!\!\!\ & & \\
{\scriptstyle E}\uparrow\phantom{a} &
\phantom{aaaaaa}\searrow\!{\scriptstyle CH}
\!\!\!\!& \\
{[[A,B]]}& \!\!\!\!\!\!\!\!\!\!\stackrel{\chi}{\longrightarrow} &
\!\!\!\! {[[SC({\bf T})\otimes A,B]]} \ ,
\end{array}
 $$
which immediately implies commutativity of (\ref{Diag1}).

To describe $\chi$ choose a sequence of continuous functions
$\kappa_n:[1,\infty)\to[0,1]$, $n\in{\bf N}$, such that
\begin{equation}\label{endpoints}
\kappa_n(1)=1,\quad \lim_{t\to\infty}\kappa_n(t)=0,\quad n\in{\bf N},
\end{equation}
and
\begin{equation}\label{continuity}
\lim_{n\to\infty}\sup_{t\in[1,\infty)}\|\kappa_{n+1}(t)-\kappa_n(t)\|=0.
\end{equation}
One way of constructing such a sequence of functions is to set
$a_n=\sum_{i=1}^n \frac{1}{i}$ and let $\kappa_n$ be the function
 $$
\kappa_n(t)=\left\lbrace\begin{array}{ll}
1\, ,&t\in[1,a_n],\\
a_n+1-t\, ,&t\in[a_n,a_n+1],\\
0\, ,&t\in[a_n+1,\infty),
\end{array}\right.
 $$
but the actual choice is not important as soon as (\ref{endpoints}) and
(\ref{continuity}) are satisfied.
Put $K(t)= \sum_{i \in {\bf N}} \kappa_i(t)$.
Denote by $P$
the projection $\sum_{i\in{\bf N}}e_{ii}$ in $l_2({\bf Z})$.
Then $PTP$ is a one-sided shift of index one.
The asymptotic homomorphism $\chi$ is then determined by the condition that
$$
\lim_{t \to \infty} \|\chi_t(f\otimes e^{2\pi ix}) - f(K(t)) PTP\| = 0 ,
$$
where $f\in C_0(0,1)$ and
$e^{2\pi ix}$ is a generator for $C({\bf T})$.

Let $\varphi=\{\varphi_t\}_{t\in[1,\infty)}:A\to B$ be an asymptotic
homomorphism. Then $CH\circ E[\varphi]$ is equivalent to
the asymptotic homomorphism
$\psi=\{\psi_t\}_{t\in[1,\infty)}:SC({\bf T})\otimes A\to B$ defined
by
 $$
\psi_t(f\otimes e^{2\pi ix}\otimes a)= T\sum_{i\in{\bf N}}
f(\kappa_i(t))\varphi_{t_i}(a)e_{ii}\, .
 $$
Define another asymptotic homomorphism
$\psi'=\{\psi'_t\}_{t\in[1,\infty)}:SC({\bf T})\otimes A\to B$ by
 $$
\psi'_t(f\otimes e^{2\pi ix}\otimes a)=PTP \sum_{i\in{\bf N}}
f(\kappa_i(t))\varphi_{t_i}(a) e_{ii}\, .
 $$
Since
 $$
\lim_{t\to\infty}\|\psi_t(f\otimes e^{2\pi ix}\otimes a)-
\psi'_t(f\otimes e^{2\pi ix}\otimes a)\|=0
 $$
for any $f\in C_0(0,1)$, $a\in A$, it follows that the asymptotic homomorphisms
$\psi$ and $\psi'$ are equivalent.
By using the freedom in the choice of the $\kappa_i$'s we can arrange that
there is a sequence $0<m_1<m_2<\ldots$ in $\bf N$ such that
 $$
\kappa_i(t)=0 \quad{\rm for\ all}\ \
t\in[t_j,t_{j+1}],\quad i=1,2,\ldots,m_j
 $$
and
 $$
\kappa_i(t)=1  \quad{\rm for\ all}\ \
t\in[t_j,t_{j+1}],\quad i\geq m_{j+1}.
 $$
Define a new sequence $s_1\leq s_2\leq s_3\leq\ldots$ in $[1,\infty)$ such
that
\begin{eqnarray*}
&s_i=0\, ,\ 0\leq i<m\, ,&\\
&s_{m_1}=s_{m_1+1}=\ldots=s_{m_2-1}=t_1\, ,&\\
&s_{m_2}=s_{m_2+1}=\ldots=s_{m_3-1}=t_2\, ,&
\end{eqnarray*}
and so on. Then $\{\varphi_{s_n}\}$ is also a discretization for
$\varphi$, so $\psi'$ is homotopic to $\psi''$, where
 $$
\psi''_t(f\otimes e^{2\pi ix}\otimes a)=PTP \sum_{i\in{\bf N}}
f(\kappa_i(t))\varphi_{s_i}(a) e_{ii}\,
 $$
asymptotically as $t \to \infty$. Since
 $$
[\chi\otimes\varphi]_t(f\otimes e^{2\pi ix}\otimes a)=PTP \sum_{i\in{\bf N}}
f(\kappa_i(t))\varphi_t(a) e_{ii}\, ,
 $$
asymptotically as $t \to \infty$, we find that
\begin{eqnarray*}
&&\lim_{t\to\infty}\sup_i\|\psi''_t(f\otimes e^{2\pi ix}\otimes
a)-[\chi\otimes\varphi]_t(f\otimes e^{2\pi ix}\otimes a)\|\\
&&\leq\lim_{t\to\infty}\|f(\kappa_i(t))\varphi_{s_i}(a)-
f(\kappa_i(t))\varphi_t(a)\|=0
\end{eqnarray*}
for any $f\in C_0(0,1)$, $a\in A$. Since elements of the form
$f\otimes e^{2\pi ix}\otimes a$ generate $SC({\bf T})\otimes A$ as a
$C^*$-algebra, it follows that
$\lim_{t\to\infty}\|\psi''_t(z)-[\chi\otimes\varphi]_t(z)\|=0$ for all
$z\in SC({\bf T})\otimes A$. Consequently the asymptotic homomorphisms $CH\circ
E[\varphi]$ and $\chi\otimes\varphi$ are homotopic.
\hfill $\Box$

\medskip

Since $\chi$ is an
isomorphism, it follows that $CH :\Ext^{-1/2} (S^2A,B) \to [[S^3A,B]]$
is surjective. But the inverse in $E$-theory of the \ah defining $\chi$
is a genuine $*$-homomorphism $\mu : SA \to S^3A \otimes M_2$ and the
naturality of the Connes--Higson construction gives us a commuting diagram
$$
\begin{array}{ccc}
\Ext^{-1/2}(SA,B)\!\!\!\!
&\stackrel{\mu^*}{\longleftarrow}&\!\!\!\!\Ext^{-1/2}(S^3A \otimes M_2,B)\\
\downarrow{\scriptstyle CH}&&\downarrow{\scriptstyle CH} \\
{[[S^2A,B]]} &  \stackrel{\mu^*}{\longleftarrow} & {[[S^4A\otimes M_2,B]]}
\end{array}
$$
We see that this proves a) of Theorem \ref{MAIN}.

\medskip

To complete the proof Theorem \ref{MAIN} it now suffices to show that
the $CH$-map of diagram (\ref{Diag1}) is injective. The rest of the paper
is devoted to this.

\section{Proof of b) of Theorem \protect\ref{MAIN}}

Given two commuting unitaries $S,T$ in a $C^*$-algebra, we define a
projection $P(S,T)$ in the $2 \times 2$ matrices over the $C^*$-algebra
generated by $S$ and $T$ in the following way. Let
$s, c_0,c_1 : [0,1] \to \bf R$ be the functions
$$
c_0(t) = |\cos (\pi t)|1_{[0,\frac{1}{2}]}(t) \ ,
\quad c_1(t) = |\cos (\pi t)|1_{(\frac{1}{2},1]}(t) \ ,
\quad s(t) = \sin (\pi t)  \ ,
$$
where $1_{[0,\frac{1}{2}]}$, $1_{(\frac{1}{2},1]}$ are the
characteristic functions of the corresponding segments.
Set $\tilde{g} = sc_0, \ \tilde{h} = sc_1$ and $\tilde{f} = s^2$.
Since $\tilde{f}, \ \tilde{g}$ and $\tilde{h}$ are continuous and
$1$-periodic they give rise to continuous functions, $f,g,h$, on $\bf T$.
Set
$$
P(S,T) \ = \ \left ( \begin{array}{cc} f(S) & g(S) + h(S)T \\
h(S)T^* + g(S) & 1 - f(S) \end{array} \right )  \ ,
$$
cf. \cite{L}. In particular, this gives us a projection
$P \in C({\bf T}^2)\otimes M_2$ when we apply the recipe to the canonical
generating unitaries of $C({\bf T}^2)$. Note that $P$ is an element of
$M_2((SC({\bf T}))^+) \subseteq M_2(C({\bf T}^2))$. In general, $P(S,T)$
is in the range of $\id_{M_2}\otimes \lambda$, where
$\lambda : (SC({\bf T}))^+ \to C^*(S,T)$ is the unital $*$-homomorphism
with
$$
\lambda( (1 - e^{2 \pi i x}) \otimes 1) \ = \ 1 - S \ , \quad
\lambda( (1 - e^{2 \pi i x}) \otimes e^{2 \pi i y}) = T - ST \ .
$$
Consider also the projection
$$
P_0 =
{\small\left ( \begin{array}{cc} 0 & 0 \\ 0 & 1  \end{array} \right )}
\in M_2 \subseteq C({\bf T}^2)\otimes M_2 \ .
$$
We can then define a map
$$
\Bott_A : \Ext^{ - 1/2} (C({\bf T}^2) \otimes A,B)_h \to
\Ext^{ - 1/2}(A,B)_h
$$
such that
$$
[\varphi] \mapsto [(\id_{M_2} \otimes \varphi) \circ b_A] - [(\id_{M_2}
\otimes \varphi) \circ b_0],
$$
where $b_A,b_0 : A \to M_2(C({\bf T}^2)) \otimes A$ are the maps $b_A(a) =
P \otimes a$ and $b_0(a) = P_0 \otimes a$, respectively. The main part of
the proof will be to establish the following.

\begin{prop}\label{mainlemma}
Let $i : SA \to C({\bf T})\otimes A$ be the
canonical embedding, $e : C({\bf T}) \otimes A \to A$ the map obtained from
evaluation at $1 \in \bf T$ and $c : A \to C({\bf T}) \otimes A$ the
$*$-homomorphism which identifies $A$ with the constant $A$-valued
functions over $\bf T$. Then
$$
- \Bott_{SA} \circ E \circ CH( [\psi] - e^* \circ c^*[\psi]) = i^*[\psi]
$$
in $\Ext^{ - 1/2}(SA,B)_h$ for every semi-invertible extension
$\psi \in \Hom(C({\bf T}) \otimes A,Q(B))$.
\end{prop}

To begin the proof of Proposition \ref{mainlemma}, observe that
$c^*([\psi] - e^* \circ c^*([\psi])) = 0$ in $\Ext^{- 1/2}(A,B)$.
We can therefore add an asymptotically split extension $\chi$ to
$c^*(\psi - e^* \circ c^*(\psi))$ such that the resulting extension
is asymptotically split. It follows that
$$
\psi^{\prime} = \psi - e^* \circ c^*(\psi ) + e^*(\chi)
$$
is a semi-invertible extension of $C({\bf T})\otimes A$ by $B$ such that
$i^*[\psi^{\prime}] = i^*[\psi]$ and $c^*(\psi^{\prime})$ is an
asymptotically split extension of $A$ by $B$. Since $CH[e^*(\chi)] =
(Se)^*(CH[\chi]) = 0$ because $\chi$ is asymptotically split, it suffices
(by using $\psi^{\prime}$ instead of $\psi$) to consider a
semi-invertible extension $\psi \in \Hom(C({\bf T}) \otimes A,Q(B))$
with the property that $c^*(\psi)$ asymptotically splits, and show
that $\Bott_{SA} \circ E \circ CH[\psi] = i^*[\psi]$. So let $\psi$
be such an extension and set $\varphi = \psi \circ i$.

\begin{lemma}\label{L3}
Let $e^{2 \pi i x}$ denote the identity function
of the circle $\bf T$. There is a unitary $U \in M(M_2(B))$ such that
$$
{\footnotesize
\left ( \begin{array}{cc} \psi(e^{2 \pi i x}f \otimes a) & {} \\
{} & 0 \end{array} \right )
}
= q_{M_2(B)}(U)
{\footnotesize
\left ( \begin{array}{cc} \psi(f \otimes a) & {} \\
{} & 0 \end{array} \right )
}
$$
for all $ f \in C({\bf T}), \ a \in A$.
\end{lemma}

{\bf Proof.} We use the well-known fact, \cite{P}, that a surjective
$*$-homomorphism between separable $C^*$-algebras admits a surjective
unital extension to a $*$-homomorphism between the multiplier algebras.
The $*$-homomorphism $\psi$ extends to a unital $*$-homomorphism
$$
\hat{\psi} : M(C({\bf T}) \otimes A) \to M(\psi(C({\bf T}) \otimes A)).
$$
Then $V = \hat{\psi}(e^{2 \pi i x} \otimes 1_A)$ is a unitary in
$M(\psi(C({\bf T}) \otimes A))$ ($1_A$ means here the unit in
$M(A)$ and hence $e^{2 \pi i x} \otimes 1_A$ is really just the
identity function of $\bf T$ considered as a unitary multiplier
of $C({\bf T}) \otimes A$). Set $D = q_B^{-1}(\psi(C({\bf T}) \otimes A))
\subseteq M(B)$. Since $q_B$ maps $D$ onto $\psi(C({\bf T}) \otimes A))$
it extends to a surjective unital $*$-homomorphism
$\widehat{q_B} : M(D) \to M(\psi(C({\bf T}) \otimes A))$.
Since ${\footnotesize\left ( \begin{array}{cc} V & {} \\
{} & V^* \end{array} \right )}$ is in the connected
component of $1$ in $M_2( M(\psi(C({\bf T}) \otimes A)))$, there is a
unitary $U \in M_2(M(D))$ such that
$$
\id_{M_2} \otimes \widehat{q_B}(U) \ = \
{\footnotesize\left ( \begin{array}{cc}  V &
{} \\ {} & V^* \end{array} \right )} \ .
$$
Note that $M(D) \subseteq M(B)$ since $B$ is an essential ideal in $D$.
We can therefore regard $U$ as a unitary in $M_2(M(B))$. It is then
clear that $U$ has the stated property.
\hfill $\Box$

\medskip

It follows from Lemma \ref{L3} that after adding $0$ to $\psi$ and
$\varphi$, we may assume that there is a unitary $w \in M(B)$ such that
\begin{equation}\label{E1}
q_B(w) \psi(f \otimes a) = \psi(e^{2 \pi i x
}f \otimes a) \ , \ f \in C({\bf T}), \ a \in A \ .
\end{equation}
Let $\{\pi_t\}_{t \in [1,\infty)} : A \to M(B)$ be an asymptotic
homomorphism such that $\psi(1 \otimes a) = q_B(\pi_t(a))$ for all
$a$ and $t$.

\begin{lemma}\label{L2}
Let $\{u_t\}_{t \in [1,\infty)}$ be a continuous
approximate unit for $B$ such that $\lim_{t\to\infty} u_t\pi_1(a) -
\pi_1(a)u_t = 0$ for all $a \in A$. There is then an increasing
continuous function $r : [1,\infty) \to [1,\infty)$ such that $r(t)
\geq t$ for all $t \in [1,\infty)$ and $\lim_{t \to \infty}
f(u_{r(t)}) \pi_1(a) - f(u_{r(t)})\pi_t(a) = 0$ for all $a \in A,
\ f \in C_0(0,1)$.
\end{lemma}

{\bf Proof.} By the Bartle--Graves selection theorem \cite{BG} there is a
continuous function $\chi : A \to M(B)$ such that $\chi(a) - \pi_1(a) \in B$
for all $A$. The same selection theorem also provides us with an
equicontinuous asymptotic homomorphism
$\pi^{\prime} = (\pi^{\prime}_t) : A \to M(B)$ such that
$\lim_{t \to \infty}\pi_t(a) - \pi^{\prime}_t(a) = 0$ for all
$ a \in A$. Let $F_1 \subseteq F_2 \subseteq F_3 \subseteq \cdots $
be a sequence of finite subsets with dense union in $A$. By using that
$\{\pi_t (a) - \chi(a) : t \in [1,n], \ a \in F_n\}$ is a compact subset
of $B$ for all $n$, it is then straightforward to construct an $r$ with
$r(t) \geq t$ such that $\lim_{t \to \infty} u_{r(t)} \pi_t(a) -
u_{r(t)}\chi(a) = 0$ for all $a \in \bigcup_n F_n$. It follows that
$\lim_{t \to \infty} u_{r(t)} \pi^{\prime}_t(a) - u_{r(t)}\chi(a) = 0$
for all $a \in \bigcup_n F_n$, and by continuity of $\chi$ and
equicontinuity of $\{\pi^{\prime}_t\}$ it follows that this actually
holds for all $a \in A$. But then $\lim_{t \to \infty} u_{r(t)}
\pi_t(a) - u_{r(t)}\pi_1(a) = 0$ since $\chi(a) - \pi_1(a) \in B$
for all $a \in A$. The fact that $\lim_{t \to \infty}
f(u_{r(t)})\pi_t(a) - f(u_{r(t)})\pi_1(a) = 0$ for all
$a \in A, \ f \in C_0(0,1)$, then follows from Weierstrass' theorem.
\hfill $\Box$

\medskip

It follows from Lemma \ref{L2} that
$CH[\psi] \in [[SC({\bf T})\otimes A, B]]$ is represented by an
asymptotic homomorphism $CH(\psi)$ such that
$$
\lim_{t \to \infty} CH(\psi)_t(f \otimes g \otimes a) - f(u_{r(t)}) g(w)
\pi_t(a) = 0
$$
for all $a \in A, \ g \in C({\bf T}), \ f \in C_0(0,1)$. By choosing
the approximate unit
$\{u_t\}$ in Lemma \ref{L2} appropriately \cite{A} we may assume that
$\lim_{t \to \infty} u_{r(t)}\pi_t(a) - \pi_t(a) u_{r(t)} = 0,
\ \lim_{t \to \infty} (1 - u_{r(t)}) (w\pi_t(a) - \pi_t(a)w) = 0$
for all $a \in A$, and $\lim_{t \to \infty} u_{r(t)} w - wu_{r(t)} = 0$.
We can therefore find a discretization $CH(\psi)_{t_i},  i \in {\bf N}$,
of $CH(\psi)$ such that

\vspace{-1.2\itemsep}
\begin{enumerate}
\item[1)] $\lim_{i \to \infty} \pi_{t_i}(a) - \pi_{t_{i+1}}(a) = 0$
for all $a \in A$,

\vspace{-1.2\itemsep}
\item[2)] $\lim_{i \to \infty} u_{r(t_i)} - u_{r(t_{i+1})} = 0$,

\vspace{-1.2\itemsep}
\item[3)] $\lim_{i \to \infty} wu_{r(t_i)} - u_{r(t_i)}w = 0 \ , $

\vspace{-1.2\itemsep}
\item[4)] $\lim_{i \to \infty} (1 - u_{r(t_i)}) (w \pi_{t_i}(a) -
\pi_{t_i}(a)w) = 0$ for all $a \in A$.
\end{enumerate}

\vspace{-1.2\itemsep}

To simplify notation, set $\pi_n = \pi_{t_n}$ and $u_n = u_{r(t_n)}$.
Set $\pi_n = \pi_{-n}$ when $n < 0$ and $\pi_0 = \pi_1$. We find that
$E \circ CH[\psi] \in \Ext^{-1/2}(SC({\bf T}^2)\otimes A,B)_h$ is
represented by a $*$-homomorphism $\underline{\Phi}$, where
$\Phi : SC({\bf T}^2) \otimes A \to {\cal L}_B(l_2({\bf Z}) \otimes B)$
is a map such that
$$
\Phi( f \otimes g \otimes h \otimes a) =
\Bigl(\sum_{n \geq 0}
f(u_n)e_{nn}\Bigr)
\Bigl(\sum_{n \in {\bf Z} } g(w) e_{nn}\Bigr) h(T)
\Bigl(\sum_{n \in {\bf Z}}
\pi_n(a) e_{nn}\Bigr)
$$
modulo ${\cal K} \otimes B$ for all $f \in C_0(0,1), \ g,h \in C({\bf T}),
\ a \in A$.

Set
$$
W = \sum_{n \in {\bf Z}} we_{nn},
\qquad U = \sum_{n \geq 0} u_n.
$$
Then $W, \ T$ and $U$ commute modulo ${\cal K} \otimes B$. Define
$\underline{\pi} : A \to Q({\cal K} \otimes B)$ such that
$\underline{\pi}(a) = \underline{\sum_n \pi_{n} (a)e_{n n}}$.
Then $\underline{\pi}$ is a $*$-homomorphism which commutes with
$\underline{U}$ and $\underline{T}$.

Let $Q \in M_2(Q(B))$ be the projection
$$
 Q = \left ( \begin{array}{cc} s^2(\underline{U}) & sc_0(\underline{U}) +
sc_1(\underline{U})\underline{W} \\ sc_1(\underline{U})\underline{W}^* +
sc_0(\underline{U}) & (c_0 + c_1)^2(\underline{U}) \end{array} \right ) \ .
$$

\begin{lemma}\label{L55}
$ - \Bott_{SA} \circ E \circ CH[\psi]$ is
represented in $\Ext^{- 1/2}(SA,B)_h$ by an extension $\lambda :
SA \to M_2(Q(B))$ such that
$$
\lambda( (1 -e^{2 \pi i x}) \otimes a) \ = \ Q
{\footnotesize\left ( \begin{array}{cc}
1 - \underline{T} & {} \\ {} & 1- \underline{T} \end{array} \right )}
{\footnotesize\left ( \begin{array}{cc} \underline{\pi}(a) & {} \\ {} &
\underline{\pi}(a) \end{array} \right )} \ ,
$$
$a \in A$.
\end{lemma}

{\bf Proof.} To simplify notation, set
$$
\tilde{U} =
{\footnotesize
\left ( \begin{array}{cc} (1 - e^{2 \pi i x})(\underline{U}) &
{} \\ {} & (1 -e^{2 \pi i x})(\underline{U}) \end{array} \right )
}.
$$
By definition $\Bott_{SA} \circ E \circ CH[\psi] = [\lambda_+] -
[\lambda_-]$, where $\lambda_{\pm} : SA \to M_2(Q(B))$ are
$*$-homomorphisms such that
$$
\lambda_+( (1 - e^{2 \pi i x}) \otimes a) \ = \ P(\underline{T},
\underline{W})\tilde{U}
{\footnotesize\left ( \begin{array}{cc} \underline{\pi}(a) &
{} \\ {} & \underline{\pi}(a) \end{array} \right )
}
$$
and
$$
\lambda_-( (1 - e^{2 \pi i x} ) \otimes a) \ = \ P_0 \tilde{U}
{\footnotesize\left ( \begin{array}{cc} \underline{\pi}(a) & {} \\ {} &
\underline{\pi}(a) \end{array} \right )} \ .
$$
Set
$$
X =  1 - (1 - P(\underline{T}, \underline{W})\tilde{U})
( 1 -P_0\tilde{U}^*) \ .
$$
Then $[\lambda_+] - [\lambda_-] = [\lambda^{\prime}]$, where
$\lambda^{\prime} : SA \to M_2(Q(B))$ is given by
$$
\lambda^{\prime}((1 - e^{2 \pi i x})\otimes a) = X
{\footnotesize\left ( \begin{array}{cc} \underline{\pi}(a) & {} \\ {} &
\underline{\pi}(a) \end{array} \right )} \ .
$$
Note that $X$ is an element in the $2\times 2$ matrices over the
$C^*$-algebra generated by $1 - \underline{T}, \ \underline{W}$ and
$(1 - e^{2 \pi i x})(\underline{U})$. In fact, if we define
$\Lambda : S \otimes C({\bf T}) \otimes S \to Q(B)$ such that
$$\Lambda ( (1 -e^{2 \pi i x} ) \otimes e^{ 2 \pi i y} \otimes
(1 - e^{2 \pi i z})) = (1 -e^{2 \pi i x})(\underline{U})\ \underline{W}\
(1 - \underline{T}),
$$
there is a quasi-unitary $d \in M_2( S \otimes C({\bf T}) \otimes S)$
such that $\id_{M_2} \otimes \Lambda(d) = X$. Here $S$ is shorthand for
the $C^*$-algebra $C_0(0,1)$. Also we remind the reader that a
quasi-unitary is an element $d$ of a $C^*$-algebra $D$ such that $1-d$ is
unitary in $D^+$. Alternatively, it is a normal element with spectrum in
$ \{ 1 - z : z \in {\bf T}\}$. Then
$$
\id_{M_2} \otimes \Lambda \otimes
{\footnotesize\left ( \begin{array}{cc}
\underline{\pi} & {} \\ {} & \underline{\pi} \end{array} \right )}
:
S \otimes C({\bf T}) \otimes S \otimes A \to M_2(Q(B))
$$
is semi-invertible, with the inverse given by the $*$-homomorphism which
results when one replaces $U$ with $\sum_{n < 0} u_{-n}$ in the definition
of $\Lambda$. Define
$$
\nu : SA \to M_2(S \otimes C({\bf T}) \otimes
S \otimes A)
$$
such that $\nu( (1 - e^{2 \pi i x}) \otimes a) =
d \otimes a$ and note that $\lambda^{\prime} = \Bigl(\id_{M_2} \otimes
\Lambda \otimes
{\footnotesize\left ( \begin{array}{cc} \underline{\pi} & {} \\
{} & \underline{\pi} \end{array} \right )}\Bigr)  \circ \nu$. Let
$\alpha$ be the automorphism of $M_2(S \otimes C({\bf T}) \otimes S
\otimes A)$ which exchanges the two suspensions by a $\pi/2$ rotation
of ${\bf R}^2$. Then
$$
- \left[\left(\id_{M_2} \otimes\Lambda \otimes
{\footnotesize\left ( \begin{array}{cc}
\underline{\pi} & {} \\ {} & \underline{\pi} \end{array} \right )
}\right)\right] =
\left[\left(\id_{M_2} \otimes \Lambda \otimes
{\footnotesize\left ( \begin{array}{cc}
\underline{\pi} & {} \\ {} & \underline{\pi} \end{array} \right )}\right)
\circ \alpha \right]
$$
in $\Ext^{ -1/2} (S \otimes C({\bf T}) \otimes S \otimes A, B)_h$.
It follows that
$$
-[\lambda^{\prime}] =  \left[\left(\id_{M_2} \otimes \Lambda \otimes
{\footnotesize\left ( \begin{array}{cc} \underline{\pi} & {} \\ {} &
\underline{\pi} \end{array} \right )}\right) \circ \alpha \circ \nu\right]
$$
in $\Ext^{ -1/2} (SA, B)_h$. Set
$$
Y =  1 - \Bigl(1 - Q
{\footnotesize\left ( \begin{array}{cc} 1 -\underline{T} & {} \\ {} &
1 -\underline{T} \end{array} \right )}\Bigr) \Bigl( 1 -P_0
{\footnotesize\left ( \begin{array}{cc} 1 -\underline{T}^* & {} \\ {} &
1 -\underline{T}^* \end{array} \right )}\Bigr) \ ,
$$
and note that $\Bigl(\id_{M_2} \otimes \Lambda \otimes
{\footnotesize\left ( \begin{array}{cc} \underline{\pi} & {} \\ {} &
\underline{\pi} \end{array} \right )}\Bigr) \circ \alpha  \circ \nu = \mu$,
where $ \mu : SA \to M_2(Q(B))$ is such that
$$
\mu((1 - e^{2 \pi i x}) \otimes a) = Y
{\footnotesize\left ( \begin{array}{cc} \underline{\pi}(a) & {} \\ {} &
\underline{\pi}(a) \end{array} \right )} \ .
$$
It follows that $[\mu] = [\lambda] - [\mu^{\prime}]$, where
$\mu^{\prime}( (1 - e^{2 \pi i x}) \otimes a) =
(1 - \underline{T})\underline{\pi}(a)$. It is easily seen that
$\mu^{\prime}$ is asymptotically split. Therefore $[\mu] = [\lambda]$.
\hfill $\Box$

\medskip

Set
$$
X = \left ( \begin{array}{cc} s(\underline{U})  & - c_0(\underline{U}) -
c_1(\underline{U})\underline{W} \\ c_0(\underline{U}) +
c_1(\underline{U})\underline{W}^* & s(\underline{U})
\end{array} \right ) \ .
$$
and
$$
Z =
\left ( \begin{array}{cc} i\underline{W_+} & 0 \\ 0 & - i\underline{W_+}
\end{array} \right )\ ,
$$
where
$$
W_+ = \sum_{n \geq 0} we_{nn}  + \sum_{n < 0} e_{nn} \in
{\cal L}_B(l_2({\bf Z})\otimes B).
$$
Then $Z$ and $X$ are unitaries in $M_2(Q(B))$. Let
$T_0 : l_2({\bf Z}) \otimes B \to l_2({\bf Z}) \otimes B$ be the unitary
$$
T_0  = \sum_{n \in {\bf Z} \backslash \{-1\}} e_{n,n+1} + we_{-1,0} \ .
$$
We can then define an extension $\lambda_1 : SA \to Q(B)$ such that
$$
\lambda_1((1 -e^{ 2 \pi i x})\otimes a) = (1 - \underline{T_0})\
\underline{\pi}(a) \ .
$$

\begin{lemma}\label{L6}
Let $\lambda : SA \to M_2(Q(B))$ be the extension
of Lemma {\rm \ref{L55}}. Then
$$
\Ad X^* \circ \lambda = \Ad Z  \circ
{\footnotesize\left (\begin{array}{cc}
\lambda_1 & {} \\ {} & 0 \end{array} \right )} \ .
$$
\end{lemma}

{\bf Proof.} Note that $\lambda$ and
${\footnotesize\left (\begin{array}{cc}
\lambda_1 & {} \\ {} & 0 \end{array} \right )}$ both extend to
unital $*$-homomorphisms $C({\bf T}) \otimes A \to M_2(Q(B))$ defined
such that they send $1 \otimes a$ to
${\footnotesize\left ( \begin{array}{cc}
\underline{\pi}(a) & {} \\ {} & \underline{\pi}(a)
\end{array} \right )}$, $ a \in A$. By considering these
extensions we see that it suffices to show that
\begin{equation}\label{F1}
X^*
{\footnotesize\left ( \begin{array}{cc} \underline{\pi}(a) & {} \\ {} &
\underline{\pi}(a) \end{array} \right )} X \ = \ Z
{\footnotesize\left ( \begin{array}{cc} \underline{\pi}(a) & {} \\ {} &
\underline{\pi}(a) \end{array} \right )} Z^*  \ ,
\end{equation}
and
\begin{eqnarray}\label{F2}
&&V^* \left( Q
{\footnotesize\left ( \begin{array}{cc} \underline{T} & {} \\
{} & \underline{T}^* \end{array} \right )} +
Q^{\perp}\right)
{\footnotesize\left ( \begin{array}{cc} \underline{\pi}(a) & {} \\
{} & \underline{\pi}(a) \end{array} \right )} V\nonumber\\
&&=
Z {\footnotesize\left ( \begin{array}{cc}  \underline{T_0} & {} \\
{} & 1 \end{array} \right )}
{\footnotesize\left ( \begin{array}{cc} \underline{\pi}(a) & {} \\
{} & \underline{\pi}(a) \end{array} \right )} Z^* \\
&&= Z \left({\footnotesize\left ( \begin{array}{cc}  1 & {} \\
{} & 0 \end{array}
\right )}
{\footnotesize\left ( \begin{array}{cc}  \underline{T_0} & {} \\ {} &
\underline{T_0}^*)  \end{array} \right)} + {\footnotesize\left (
\begin{array}{cc}  0 & {} \\ {} & 1 \end{array} \right )} \right)
{\footnotesize\left (
\begin{array}{cc} \underline{\pi}(a) & {} \\ {} & \underline{\pi}(a)
\end{array} \right )} Z^* .\nonumber
\end{eqnarray}
To simplify the verification, observe that $W_+T_0 = TW_+$ from
which it follows that $Z{\footnotesize\left ( \begin{array}{cc}
\underline{T_0} & {} \\ {} & \underline{T_0}^*
\end{array} \right )}Z^* =
{\footnotesize\left ( \begin{array}{cc}  \underline{T} &
{} \\ {} & \underline{T}^* \end{array} \right )}$. Since $X$ clearly
commutes with
${\footnotesize\left ( \begin{array}{cc}  \underline{T} & {} \\ {} &
\underline{T}^* \end{array} \right )}$ and $Z$ with
${\footnotesize\left ( \begin{array}{cc} 1 & {} \\
{} & 0 \end{array} \right )}$
we see that (\ref{F2}) will follow from $(\ref{F1})$ and
\begin{equation}\label{F3}
X^*Q  X = {\footnotesize\left ( \begin{array}{cc} 1 & {} \\ {} &
0 \end{array} \right )} \ .
\end{equation}
Thus we need only check (\ref{F1}) and (\ref{F3}), and we leave
that to the reader.
\hfill $\Box$

\medskip

All in all we now have that $- \Bott_{SA} \circ E \circ CH[\psi] =
[\lambda_1]$ in $\Ext^{ - 1/2} (SA,B)_h$. Define
$$\kappa :
SA \to Q(B)
\quad {\rm by}\quad
\kappa((1- e^{2 \pi i x}) \otimes a) =
(1 - \underline{T})\, \underline{\pi}(a).
$$
The extension
$\kappa$ is asymptotically
split and hence $[\lambda_1] = [\lambda_1] - [\kappa]$. Since
$[\lambda_1] - [\kappa] = [\mu]$, where $\mu : SA \to Q(B)$ is given by
$\mu((1 -  e^{2 \pi i x}) \otimes a) = (1 - T_0T^*) \underline{\pi}(a)$
and since $T_0T^* = \sum_{n \neq -1} e_{nn} + we_{-1,-1}$, we see that
$[\lambda_1] = [\varphi]$. This completes the proof of Proposition
\ref{mainlemma}.

\begin{cor}\label{C1}
The map
$$
CH : \Ext^{ -1/2}(SA,B)_h \to [[S^2A,B]]
$$
is injective on $i^*(\Ext^{ - 1/2}(C({\bf T})\otimes A,B)_h)$.
\end{cor}

{\bf Proof.} Let $\psi \in \Ext^{ -1/2}(C({\bf T}) \otimes A, B)$ and
assume that $CH(i^*[\psi]) = 0$. By the naturality of the Connes--Higson
construction this implies that
$$
(Si)^*(CH[\psi] - (Se)^* \circ (Sc)^*(CH[\psi])) = CH(i^*(\psi)) = 0
$$
in $[[S^2A,B]]$. But the split exactness of the functor
$[[S - , B]]$, \cite{DL}, implies then that
$$
0 = CH[\psi] - (Se)^* \circ
(Sc)^*(CH[\psi]) = CH([\psi] - e^* \circ c^*[\psi]).
$$
And then $i^*[\psi] = 0$ by Proposition \ref{mainlemma}.
\hfill $\Box$

\medskip

\begin{lemma}\label{L4}
The map
$$
(Si)^* : \Ext^{- 1/2}(SC({\bf T}) \otimes A,B)_h
\to  \Ext^{- 1/2}(S^2A,B)_h
$$
is surjective.
\end{lemma}

{\bf Proof.} To prove this we shall identify $S^2 = C_0({\bf R}^2)$ with
$C_0(D)$, where $D \ = \ {\bf R}^2 \backslash \{(0,y) \in {\bf R}^2 \ :
y \geq 0\}$ and $SC({\bf T})$ with $C_0({\bf R}^2 \backslash \{0\})$.
It is easy to see that there is a continuous map
$F : [0,1] \times {\bf R}^2 \to {\bf R}^2$ such that $F(0, -)$ is a
homeomorphism $\mu = F(0,-) : {\bf R}^2 \to D$; \ $F(1,z) = z, \ z \in
{\bf R}^2$, and $F^{-1}(K)$ is compact for every compact subset $K$ of $D$.
It follows that $f \mapsto f \circ \mu^{-1}$ is an endomorphism of $C_0(D)$
which is homotopic to $\id_{C_0(D)}$. Hence if
$\varphi \in \Hom(S^2A, Q(B))$ is a semi-invertible extension,
$[\varphi] = [\chi]$ in $\Ext^{ - 1/2}(S^2A,B)$, where $\chi(f) =
\varphi(f \circ \mu^{-1})$. Define $\psi : SC({\bf T}) \to Q(B)$ by
$\psi(g) = \varphi(g \circ \mu^{-1})$. Then $(Si)^*[\psi] = [\varphi]$.
\hfill $\Box$

\medskip

Lemma \ref{L4} and Corollary \ref{C1} in combination prove that the
$CH$-map of diagram (\ref{Diag1}) is injective. This completes the proof
of b) of Theorem \ref{MAIN}.
\hfill $\Box$

\vspace{2cm}
\parbox{7cm}{V. M. Manuilov\\
Dept. of Mech. and Math.,\\
Moscow State University,\\
Moscow, 119899, Russia\\
e-mail: manuilov@mech.math.msu.su
}
\hfill
\parbox{6cm}{K. Thomsen\\
Institut for matematiske fag,\\
Ny Munkegade, 8000 Aarhus C,\\
Denmark\\
e-mail: matkt@imf.au.dk
}

\end{document}